\documentclass[12pt]{article}%
\pdfoutput=1
\usepackage{amssymb}
\usepackage{amsfonts}
\usepackage{amsmath}
\usepackage{graphicx}
\usepackage{amscd}%
\providecommand{\U}[1]{\protect\rule{.1in}{.1in}}
\newtheorem{theorem}{Theorem}[section]
\newtheorem{proposition}[theorem]{Proposition}
\newtheorem{lemma}[theorem]{Lemma}
\newtheorem{corollary}[theorem]{Corollary}

\begin{document}

\title{A uniform reconstruction formula in integral geometry}
\author{V. P. Palamodov\\Tel Aviv University}
\date{November 2011}
\maketitle

\textbf{Abstract: }A general method for analytic inversion in integral
geometry is proposed. All classical and some new reconstruction formulas of
Radon-John type are obtained by this method. No harmonic analysis and no PDE
is used.



\section{Introduction}

We present a uniform reconstruction method for a class of Funk-Radon type
integral transforms extending results of \cite{P2} for arbitrary dimensions.
The reconstruction does not include summation of an infinite series and looks
as the standard (John type) inversion for the Radon transform. We specify this
method for classical and new regular acquisition geometries. The condition of
regularity (see the next Section) is necessary for\ an inversion operator to
be bounded in a Sobolev space scale, but it is not sufficient. We shall see
that the existence of a elementary exact reconstruction formula depends on
vanishing of some integrals of\ rational forms on an algebraic manifold. In
the last Section, we discuss reconstruction in Euclidean space from data of
integrals over a family of spheres. This subject is in focus of recent
research, see papers \cite{F2},\cite{K1},\cite{K2},\cite{XW},\cite{N} and
references therein.

\section{Geometry and integrals}

Let $X$ and $\Sigma$ be smooth $n$ dimensional manifolds $(n>1),$ let $Z$ be a
smooth closed hypersurface in $X\times\Sigma$ and $p:Z\rightarrow
X,\ \pi:Z\rightarrow\Sigma$ be natural projections. We suppose that there
exists a real smooth function $\Phi$ in $X\times\Sigma$ (called generating
function)\ such that $Z=\left\{  \left(  x,\sigma\right)  ;\ \Phi\left(
x,\sigma\right)  =0\right\}  $ and $\mathrm{d}_{x}\Phi\neq0$ on $Z$. Suppose that

(\textbf{i})\ \textit{The map }$\pi$ \textit{has rank }$n$\textit{\ and the
mapping }$P:N^{\ast}\left(  Z\right)  \rightarrow T^{\ast}\left(  X\right)
$\textit{\ is a local diffeomorphism}. Here, $N^{\ast}\left(  Z\right)  $
denotes the conormal bundle of $Z$ and $P\left(  x,\sigma;\nu_{x},\nu_{\sigma
}\right)  =\left(  x,\nu_{x}\right)  \in T^{\ast}\left(  X\right)  .$ It
follows that a set $Z\left(  \sigma\right)  =\pi^{-1}\left(  \sigma\right)
=\left\{  x;\ \Phi\left(  x,\sigma\right)  =0\right\}  $ is a smooth
hypersurface in $X,$ and for any point $x\in X$ and for any tangent hyperplane
$h\subset T_{x}\left(  X\right)  $ there is a locally unique hypersurface
$Z\left(  \lambda,\omega\right)  $ through $x$ tangent to $h.$

\begin{proposition}
\label{p1}For an arbitrary generating function $\Phi$ property (\textbf{i}) is
equivalent to the condition: $\mathrm{\det}\left(  \mathrm{d}_{x,t}%
\mathrm{d}_{\sigma,\tau}\Psi\right)  \neq0\ $where $\Psi\left(  x,t;\sigma
,\tau\right)  =t\tau\Phi\left(  x,\sigma\right)  ,$ $t,\tau\in\mathbb{R}%
,t\tau>0$ for any local coordinate system $x_{1},...,x_{n}$ in $X\ $and any
local coordinate system $\sigma_{1},...,\sigma_{n}$ in $\Sigma.$
\end{proposition}

For a proof see \cite{P1}, Proposition 1.1.

\textbf{Definition.} We call a generating function $\Phi$ \textit{regular} if
it satisfies conditions (\textbf{i}) and

(\textbf{ii}) there are no conjugated points, that is the equations
$\Phi\left(  x,\sigma\right)  =\Phi\left(  y,\sigma\right)  $ and
\textrm{d}$_{\sigma}\Phi\left(  x,\sigma\right)  =\mathrm{d}_{\sigma}%
\Phi\left(  y,\sigma\right)  $ are fulfilled for no $x\neq y\in X,$ $\sigma
\in\Sigma$.

Suppose that $X$ is oriented and $\mathrm{g}$ is a Riemannian metric in $X;$
let $\mathrm{d}V$ be the oriented Riemannian volume form. Consider the
integral%
\[
M_{\Phi}f\left(  \sigma\right)  =\int\delta\left(  \Phi\left(  x,\sigma
\right)  \right)  f\mathrm{d}S=\int_{Z\left(  \sigma\right)  }\frac
{f\mathrm{d}V}{\mathrm{d}\Phi\left(  x,\sigma\right)  }=\int_{Z\left(
\sigma\right)  }fq
\]
for an arbitrary continuous function $f$ compactly supported\ in $X.$ The
quotient $f\mathrm{d}V/\mathrm{d}\Phi$ denotes an arbitrary $n-1$ form $q$
such that $\mathrm{d}\Phi\wedge q=f\mathrm{d}V.$ It is defined up to a term
$h\mathrm{d}\theta$ where $h$ is a continuous function. An orientation of a
hypersurface $Z\left(  \sigma\right)  $ is defined by means of the form
$\mathrm{d}\theta$ and the integral of $q$ over $Z\left(  \sigma\right)  $ is
uniquely defined.$\ $We call the operator $M_{\Phi}$ \textit{Minkowski-Funk}
transform generated by $\Phi.$ This transform can be written in terms of
Riemannian integral as follows:
\begin{equation}
M_{\Phi}f\left(  \sigma\right)  =\int_{Z\left(  \sigma\right)  }%
\frac{f\mathrm{d}_{\mathrm{g}}S}{\left\vert \nabla_{\mathrm{g}}\Phi\left(
x,\sigma\right)  \right\vert } \label{1}%
\end{equation}
Here, $\mathrm{d}_{\mathrm{g}}S$ is the Riemannian $n-1$ surface element and
$\left\vert \nabla_{\mathrm{g}}a\right\vert \doteq\sqrt{\mathrm{g}\left(
\mathrm{d}a\right)  }$ is the Riemannian gradient of a function $a$. The
function $\left\vert \nabla_{\mathrm{g}}\Phi\right\vert $ vanishes nowhere
since of (\textbf{i}). Suppose that the gradient factorizes through $X$ and
$\Sigma$ that is $\left\vert \nabla_{\mathrm{g}}\Phi\left(  x,\sigma\right)
\right\vert =m\left(  x\right)  \mu\left(  \sigma\right)  $ for some positive
continuous functions $m$ in $X$ and $\mu$ in $\Sigma.\ $Then data of the
Minkowski-Funk transform is equivalent to data of Riemann hypersurface
integrals%
\[
R_{\mathrm{g}}f\left(  \sigma\right)  =\int_{Z\left(  \sigma\right)
}f\mathrm{d}_{\mathrm{g}}S,\sigma\in\Sigma
\]
since $R_{\mathrm{g}}f\left(  \sigma\right)  =\mu\left(  \sigma\right)
M_{\Phi}\left(  mf\right)  \left(  \sigma\right)  $. The reconstruction
problem of a function $f$ from Riemann integrals $R_{\mathrm{g}}f$ is then
reduced to inversion of the operator $M_{\Phi}.$

We say that a generating function $\Phi$ is \textit{resolved} if
$\Sigma=\mathbb{R}\times\mathbb{\mathrm{S}}^{n-1},$ and $\Phi$ has the form
$\Phi\left(  x;\lambda,\omega\right)  =\theta\left(  x,\omega\right)
-\lambda,\ \lambda\in\mathbb{R},\ \omega\in\mathbb{\mathrm{S}}^{n-1}$ for a
smooth function $\theta$ on $X\times\mathrm{S}^{n-1}\ $where $\mathrm{S}%
^{n-1}$ denotes the unit sphere in an $n$ dimensional space. Note that the map
$p:Z\rightarrow X$ is always proper for a resolved generating function. This
property guarantees that the functions $M_{\Phi}f$ and$\ R_{\mathrm{g}}f$ have
compact support in $\Sigma.$ The operator $M_{\Phi}$ fulfils the range
conditions similar to that of the Radon transform.

\begin{proposition}
Let $\Phi=\theta-\lambda$ be a resolved regular generating function and
$\theta\left(  x,\omega\right)  $ be a polynomial function of $\omega$ of
order $m.$ Then for an arbitrary integrable function $f$ in $X$ with compact
support and for an arbitrary polynomial $p\left(  \lambda\right)  $ of order
$k$ the integral%
\[
\int p\left(  \lambda\right)  M_{\Phi}f\left(  \lambda,\omega\right)
\mathrm{d}\lambda
\]
is a polynomial of $\omega$ of order $\leq mk.$
\end{proposition}

\textbf{Proof. }We have
\[
\int p\left(  \lambda\right)  M_{\Phi}f\left(  \lambda,\omega\right)
\mathrm{d}\lambda=\int p\left(  \lambda\right)  \int_{\theta=\lambda}%
\frac{f\mathrm{d}V}{\mathrm{d}\theta}\mathrm{d}\lambda=\int_{X}p\left(
\theta\left(  x,\omega\right)  \right)  f\left(  x\right)  \mathrm{d}V
\]
where $p\left(  \theta\left(  x,\omega\right)  \right)  $ is a polynomial of
$\omega$ of order $\leq mk.\ \blacktriangleright$

\section{Main theorem}

For a real smooth function $f$ in a manifold $X$ and a natural $n$ we consider
singular integrals%
\[
I_{n\pm}\left(  \rho\right)  =\int_{X}\frac{\rho}{\left(  f\pm i0\right)
^{n}}=\lim_{\varepsilon\searrow0}\int_{X}\frac{\rho}{\left(  f\pm
i\varepsilon\right)  ^{n}},
\]
for a smooth density $\rho$ with compact support. If $\mathrm{d}f\neq0$ on the
zero set of $f,$ then these limits exist and the functionals $I_{n\pm}$ are
generalized functions in $X.$ The functional
\[
\left(  P\right)  \int_{X}\frac{\rho}{f^{n}}\doteq\operatorname{Re}%
I_{n+}\left(  \rho\right)  =\operatorname{Re}I_{n-}\left(  \bar{\rho}\right)
\]
is called a\textit{\ principal value} integral. For a resolved regular
generating function $\Phi=\theta-\lambda$ we define the function on $X\times
X\backslash\{\mathrm{diag}\}$%
\[
\Theta_{n}\left(  x,y\right)  =\int_{\mathbb{\mathrm{S}}^{n-1}}\frac
{\mathrm{d}\omega}{\left(  \theta\left(  x,\omega\right)  -\theta\left(
y,\omega\right)  -i0\right)  ^{n}}%
\]
where $\mathrm{d}\omega$ is the Euclidean volume form on $\mathbb{\mathrm{S}%
}^{n-1}$. The singular integral converges since by (\textbf{ii})
the$\ \mathrm{d}_{\omega}\left(  \theta\left(  x,\omega\right)  -\theta\left(
y,\omega\right)  \right)  \neq0$ as $\theta\left(  x,\omega\right)
-\theta\left(  y,\omega\right)  =0.$

\begin{theorem}
\label{NR}Let $\Phi=\theta-\lambda$ be a regular resolved generating function
in $X\times\Sigma$ and $f\in L_{2}\left(  X\right)  $ be an arbitrary function
with compact support. If $n$ is even and $\operatorname{Re}\Theta_{n}\left(
x,y\right)  =0$ for any $x\neq y\in X,$ a reconstruction from data of
$M_{\Phi}f$ is given by the formula:%
\begin{align}
f\left(  x\right)   &  =-\frac{\left(  n-1\right)  !}{\left(  2\pi i\right)
^{n}}\frac{1}{D_{n}\left(  x\right)  }\left(  P\right)  \int_{\Sigma}%
\frac{M_{\Phi}f\left(  \lambda,\omega\right)  \mathrm{d}\lambda\mathrm{d}%
\omega}{\left(  \theta\left(  x,\omega\right)  -\lambda\right)  ^{n}}%
\label{9}\\
&  =-\frac{1}{\left(  2\pi i\right)  ^{n}D_{n}\left(  x\right)  }%
\int_{\mathrm{S}^{n-1}}\left(  P\right)  \int_{\mathbb{R}}\frac{\partial
^{n-1}}{\partial\lambda^{n-1}}M_{\Phi}f\left(  \lambda,\omega\right)
\frac{\mathrm{d}\lambda\mathrm{d}\omega}{\theta\left(  x,\omega\right)
-\lambda}\nonumber
\end{align}
If $n$ is odd and $\operatorname{Im}\Theta_{n}\left(  x,y\right)  =0$ for
$x\neq y,$ the function can be reconstructed by%
\begin{align}
f\left(  x\right)   &  =\frac{1}{2\left(  2\pi i\right)  ^{n-1}D_{n}\left(
x\right)  }\int_{\Sigma}\delta^{\left(  n-1\right)  }\left(  \theta\left(
x,\omega\right)  -\lambda\right)  M_{\Phi}f\left(  \lambda,\omega\right)
\mathrm{d}\lambda\mathrm{d}\omega\label{2}\\
&  =\frac{1}{2\left(  2\pi i\right)  ^{n-1}D_{n}\left(  x\right)  }%
\int_{\mathrm{S}^{n-1}}\left.  \frac{\partial^{n-1}}{\partial\lambda^{n-1}%
}M_{\Phi}f\left(  \lambda,\omega\right)  \right\vert _{\lambda=\theta\left(
x,\omega\right)  }\mathrm{d}\omega\nonumber
\end{align}
where%
\[
D_{n}\left(  x\right)  =\frac{1}{\left\vert \mathrm{S}^{n-1}\right\vert }%
\int_{\mathrm{S}^{n-1}}\frac{\mathrm{d}\omega}{\left\vert \nabla_{\mathrm{g}%
}\theta\left(  x,\omega\right)  \right\vert ^{n}}%
\]
The integrals (\ref{9}) and (\ref{2}) converge in mean on any compact set in
$X$.
\end{theorem}

\textbf{Remark 1.} A more invariant form of (\ref{9}) or (\ref{2}) is the
reconstruction of the form $f\mathrm{d}V:$%
\[
f\mathrm{d}V=-\frac{1}{\left(  2\pi i\right)  ^{n}}\frac{\mathrm{d}V}{D_{n}%
}\int...\mathrm{d}\omega
\]
The quotient $\mathrm{d}V/D_{n}$ depends only on the conformal class of the
Riemannian metric $\mathrm{g}.$

\textbf{Remark 2. }The condition (\textbf{ii}) is not required when $n=2$.

\textbf{Remark 3. }G. Beylkin studied "the generalized Radon transform"
\cite{B}, which coincides with the operator $M_{\Phi}$ in the Euclidean space.
He constructed a Fourier integral operator parametrix for this operator and
reduced inversion of this operator to solution of a Fredholm equation. Two
dimensional case of Theorem \ref{NR}\ was stated in \cite{P2}.

\textbf{Remark 4. }When the condition on $\Theta_{n}$ fails, the right-hand
side $R_{n}\ $of (\ref{9}), respectively (\ref{2}), gives anyway a
reconstruction up to a compact operator. More precisely, let $K$ be an
arbitrary compact set in $X$ and $I_{K}$ be the indicator function of $K.$
Then the equation holds $I_{K}Rf=f+C_{K}f$ for any function $f$ with support
in $K$ where $C_{K}$ is a compact operator in $L_{2}\left(  K\right)  .$

\textbf{Proof of Theorem}.\textbf{\ }

\begin{lemma}
\label{CO}The integral transform
\begin{align*}
I_{n}f\left(  x\right)   &  =\left(  P\right)  \int_{\mathrm{S}^{n-1}}%
\int_{\mathbb{R}}\frac{M_{\Phi}f\left(  \lambda,\varphi\right)  }{\Phi
^{2}\left(  x;\lambda,\varphi\right)  }\mathrm{d}\lambda\mathrm{d}%
\varphi,\ \text{for even}\ n\\
I_{n}f\left(  x\right)   &  =\int_{\mathrm{S}^{n-1}}\left.  \frac
{\partial^{n-1}}{\partial\lambda^{n-1}}Mf\left(  \lambda,\omega\right)
\right\vert _{\lambda=\theta\left(  x,\omega\right)  }\mathrm{d}%
\omega,\ \text{for odd}\ n
\end{align*}
is a continuous operator $L_{2}\left(  X\right)  _{\mathrm{comp}}\rightarrow
L_{2}\left(  X\right)  _{\mathrm{loc}}.$
\end{lemma}

A proof of the Lemma is given in \cite{P1}.\ For even $n$ and an arbitrary
$x\in X$ and a\ function $f$ that vanishes in a neighborhood of $x$ we
calculate%
\begin{align*}
I_{n}f\left(  x\right)   &  \doteq\left(  P\right)  \int_{\Sigma}\frac
{M_{\Phi}f\left(  \lambda,\omega\right)  \mathrm{d}\lambda\mathrm{d}\omega
}{\Phi^{n}\left(  x;\lambda,\omega\right)  }=\int_{\mathrm{S}^{n-1}}%
\mathrm{d}\omega\left(  P\right)  \int_{\mathbb{R}}\int_{\theta\left(
y,\omega\right)  =\lambda}f\left(  y\right)  q\frac{\mathrm{d}\lambda}{\left(
\theta\left(  x,\omega\right)  -\lambda\right)  ^{n}}\\
&  =\int_{X}\left(  \left(  P\right)  \int_{\mathrm{S}^{n-1}}\frac
{\mathrm{d}\omega}{\left(  \theta\left(  x,\omega\right)  -\theta\left(
y,\omega\right)  \right)  ^{n}}\right)  f\mathrm{d}\theta\wedge q=\int
_{X}\operatorname{Re}\Theta_{n}\left(  x,y\right)  f\left(  y\right)
\mathrm{d}V\left(  y\right)
\end{align*}
Here, the relation $\mathrm{d}\lambda=\mathrm{d}\theta$ holds in $Z$ and the
equation \textrm{d}$\theta\wedge q=\mathrm{d}V$ is fulfilled in $X$ by
definition. Thus the function $\Theta_{n}$ is the off-diagonal kernel of the
operator $I_{n}.$ It vanishes since of the assumption. Therefore $\Theta
_{n}\left(  x,y\right)  $ is supported in the diagonal and according to Lemma
\ref{CO} we have $\Theta_{n}\left(  x,y\right)  =a_{n}\left(  x\right)
\delta_{x}\left(  y\right)  $ for a locally bounded function $a_{n}$ in $X.$

If $n$ is odd we have%
\begin{align*}
I_{n}f\left(  x\right)   &  \doteq\int_{\mathrm{S}^{n-1}}\int_{\mathbb{R}%
}\delta^{\left(  n-1\right)  }\left(  \theta\left(  x,\omega\right)
-\lambda\right)  M_{\Phi}f\left(  \lambda,\omega\right)  \mathrm{d}%
\lambda\mathrm{d}\omega\\
&  =\int_{\mathbb{\mathrm{S}}}\mathrm{d}\omega\frac{\partial^{n-1}}%
{\partial\lambda^{n-1}}\left(  \int_{Z\left(  \lambda,\omega\right)
}fq\right)  _{\lambda=\theta\left(  x,\omega\right)  }=\int\Theta_{n}\left(
x,y\right)  f\left(  y\right)  \mathrm{d}V\left(  y\right) \\
\int_{Z\left(  \lambda,\omega\right)  }f\left(  y\right)  q  &  =\frac{1}{\pi
}\operatorname{Im}\int\frac{f\left(  y\right)  \mathrm{d}V}{\lambda
-\theta\left(  y,\omega\right)  -i0}%
\end{align*}%
\begin{align*}
\int\left.  \frac{\partial^{n-1}}{\partial\lambda^{n-1}}\left(  \int_{Z\left(
\lambda,\omega\right)  }fq\right)  \right\vert _{\lambda=\theta\left(
x,\omega\right)  }\mathrm{d}V  &  =\frac{\left(  n-1\right)  !}{\pi
}\operatorname{Im}\int\left[  \frac{1}{\left(  \theta\left(  x,\omega\right)
-\theta\left(  y,\omega\right)  -i0\right)  ^{n}}\right]  \mathrm{d}%
\omega~f\mathrm{d}V\\
&  =\frac{\left(  n-1\right)  !}{\pi}\int\operatorname{Im}\Theta_{n}\left(
x,y\right)  f\mathrm{d}V
\end{align*}
that is $\pi^{-1}\left(  n-1\right)  !\operatorname{Im}\Theta_{n}$ is the
kernel of $I_{n}$ which vanishes off the diagonal by the assumption. By Lemma
\ref{CO} we again conclude that $\pi^{-1}\left(  n-1\right)
!\operatorname{Im}\Theta_{n}=a_{n}\left(  x\right)  \delta_{x}\left(
y\right)  $ for a locally bounded function $a_{n}.$

Next we calculate the function $a_{n}.$ Choose a smooth function $e_{0}$ of
one variable with support in $\left[  -1,1\right]  $ such that $e_{0}\left(
0\right)  =1$ and set $e_{\varepsilon}\left(  x\right)  =e_{0}\left(
\left\vert x\right\vert ^{2}/\varepsilon^{2}\right)  $ for $x\in\mathbb{R}%
^{n}$ and any $\varepsilon>0.$ Take a point $x_{0}\in X$ and show that
\begin{equation}
\operatorname{Re}\int_{X}\mathrm{d}V\int_{\mathbb{\mathrm{S}}^{n-1}}%
\frac{e_{\varepsilon}\left(  x-x_{0}\right)  \mathrm{d}\omega}{\left(
\theta\left(  x,\omega\right)  -\theta\left(  x_{0},\omega\right)  -i0\right)
^{n}}\rightarrow a_{n}\left(  x_{0}\right) \nonumber
\end{equation}
for $n$ even and%
\begin{equation}
\frac{\left(  n-1\right)  !}{\pi}\operatorname{Im}\int_{X}\mathrm{d}%
V\int_{\mathbb{\mathrm{S}}^{n-1}}\frac{e_{\varepsilon}\left(  x-x_{0}\right)
\mathrm{d}\omega}{\left(  \theta\left(  x,\omega\right)  -\theta\left(
x_{0},\omega\right)  -i0\right)  ^{n}}\rightarrow a_{n}\left(  x_{0}\right)
\nonumber
\end{equation}
for $n$ odd as $\varepsilon\rightarrow0.$ We can change order of integrals and
integrate first over $X.$

\begin{lemma}
\label{D}If $n$ is even, we have for any$\ x_{0}\in X,$ arbitrary $\omega
\in\mathrm{S}^{n-1}$ and small $\varepsilon$
\begin{equation}
a_{n}\left(  x_{0},\omega\right)  \doteq\operatorname{Re}\int_{X}%
\frac{e_{\varepsilon}\left(  x-x_{0}\right)  \mathrm{d}V}{\left(
\theta\left(  x,\omega\right)  -\theta\left(  x_{0},\omega\right)  -i0\right)
^{n}}=\frac{\left(  -1\right)  ^{n/2-1}\pi^{\left(  n+1\right)  /2}}%
{\Gamma\left(  \left(  n+1\right)  /2\right)  }\frac{1}{\left\vert
\nabla\theta\left(  x_{0},\omega\right)  \right\vert ^{n}}+o\left(  1\right)
\label{19}%
\end{equation}
where $o\left(  1\right)  \leq C\varepsilon^{1/2}\log1/\varepsilon$ where $C$
does not depend on $\omega.$ For odd $n$ we have%
\[
a_{n}\left(  x_{0},\omega\right)  \doteq\frac{\left(  n-1\right)  !}{\pi
}\operatorname{Im}\int_{X}\frac{e_{\varepsilon}\left(  x-x_{0}\right)
\mathrm{d}V}{\left(  \theta\left(  x,\omega\right)  -\theta\left(
x_{0},\omega\right)  -i0\right)  ^{n}}=\frac{\left(  -1\right)  ^{m}\left(
2\pi\right)  ^{n}}{\left(  n-1\right)  !}\frac{1}{\left\vert \nabla
\theta\left(  x_{0},\omega\right)  \right\vert ^{n}}+O\left(  \varepsilon
\right)
\]

\end{lemma}

Taking the limit and integrating (\ref{19}) over $\mathrm{S}^{n-1}$ yields for
even $n$ the equation
\begin{align*}
a_{n}\left(  x_{0}\right)   &  =\lim_{\varepsilon\rightarrow0}\int
a_{n}\left(  x_{0},\omega\right)  \mathrm{d}\omega=\left(  -1\right)
^{n/2-1}\frac{\pi^{\left(  n+1\right)  /2}}{\Gamma\left(  \left(  n+1\right)
/2\right)  }\int\frac{\mathrm{d}\omega}{\left\vert \nabla\theta\left(
x_{0},\omega\right)  \right\vert ^{n}}\\
&  =-\frac{\left(  2\pi i\right)  ^{n}}{\left(  n-1\right)  !}\frac
{1}{\left\vert \mathrm{S}^{n-1}\right\vert }\int\frac{\mathrm{d}\omega
}{\left\vert \nabla\theta\left(  x_{0},\omega\right)  \right\vert ^{n}}%
=-\frac{\left(  2\pi i\right)  ^{n}}{\left(  n-1\right)  !}D_{n}\left(
x_{0}\right)
\end{align*}
which implies (\ref{9}). For odd $n$ we obtain%
\[
a_{n}\left(  x_{0}\right)  =\lim_{\varepsilon\rightarrow0}\int a_{n}\left(
x_{0},\omega\right)  \mathrm{d}\omega=2\left(  2\pi i\right)  ^{n-1}\frac
{1}{\left\vert \mathrm{S}^{n-1}\right\vert }\int\frac{\mathrm{d}\omega
}{\left\vert \nabla_{\mathrm{g}}\theta\left(  x_{0},\omega\right)  \right\vert
^{n}}=2\left(  2\pi i\right)  ^{n-1}D_{n}\left(  x_{0}\right)
\]
This yields (\ref{2}). This completes the proof of Theorem \ref{NR}%
.$\ \blacktriangleright$

\textbf{Proof of Lemma. }We show\ first that the $\theta$ can be replaced by a
linear function. Assume for simplicity that $x_{0}=0$ and $\partial
\theta\left(  x_{0},\omega\right)  /\partial x_{1}=\left\vert \nabla
_{\mathrm{g}}\theta\right\vert ,\partial\theta\left(  x_{0},\omega\right)
/\partial x_{2}=...=\partial\theta\left(  x_{0},\omega\right)  /\partial
x_{n}=0$ for a coordinate system $x_{1},...,x_{n}$ in a neighborhood of
$x_{0}$ such that $\mathrm{g}\left(  x_{0};\partial/\partial x_{i}%
,\partial/\partial x_{j}\right)  =\delta_{ij}.$ We have then $\mathrm{d}%
V=v\mathrm{d}x,\ \mathrm{d}x=\mathrm{d}x_{1}...\mathrm{dx}_{n}$ where $v$ is a
smooth function such that $v\left(  x_{0}\right)  =1.$

If\textbf{\ }$n$ is even, we integrate by parts $n$ times with respect to
$x_{1}:$%
\begin{align}
\operatorname{Re}\int\frac{e_{\varepsilon}\left(  x-x_{0}\right)  \mathrm{d}%
V}{\left(  \theta\left(  x,\omega\right)  -\theta\left(  x_{0},\omega\right)
-i0\right)  ^{n}}  &  =\frac{1}{\left(  n-1\right)  }\operatorname{Re}%
\int\frac{d_{1}\left(  x\right)  \mathrm{d}x}{\left(  \theta\left(
x,\omega\right)  -\theta\left(  x_{0},\omega\right)  -i0\right)  ^{n-1}%
}\label{22}\\
...  &  =\frac{1}{\left(  n-1\right)  !}\int\log\left\vert \theta\left(
x,\omega\right)  -\theta\left(  x_{0},\omega\right)  \right\vert d_{n}\left(
x\right)  \mathrm{d}x\nonumber
\end{align}%
\begin{align*}
d_{1}  &  =\frac{\partial}{\partial x_{1}}\frac{e_{\varepsilon}\left(
x-x_{0}\right)  v\left(  x\right)  }{\left\vert \nabla_{\mathrm{g}}%
\theta\left(  x\right)  \right\vert }=\frac{\partial e_{\varepsilon}}{\partial
x_{1}}\frac{v}{\left\vert \nabla_{\mathrm{g}}\theta\left(  x,\omega\right)
\right\vert }-\frac{e_{\varepsilon}\left(  \partial_{1}\left\vert
\nabla_{\mathrm{g}}\theta\left(  x,\omega\right)  \right\vert -v\partial
_{1}\left\vert \nabla_{\mathrm{g}}\theta\left(  x,\omega\right)  \right\vert
\right)  }{\left\vert \nabla_{\mathrm{g}}\theta\left(  x,\omega\right)
\right\vert ^{2}}\\
&  ...\\
d_{n}  &  =\frac{\partial}{\partial x_{1}}\frac{d_{n-1}\left(  x\right)
}{\left\vert \nabla_{\mathrm{g}}\theta\left(  x\right)  \right\vert }%
=\frac{\partial^{n}e_{\varepsilon}}{\partial x_{1}^{n}}\frac{v}{\left\vert
\nabla_{\mathrm{g}}\theta\left(  x,\omega\right)  \right\vert ^{n}}+...
\end{align*}
where omitted terms only include derivatives of $e_{\varepsilon}$ of order
$<n.$ Changing the variables $x=\varepsilon y$ we get%
\begin{equation}
d_{n}\left(  y\right)  =\varepsilon^{-n}\frac{\partial^{n}e}{\partial
y_{1}^{n}}\frac{v\left(  \varepsilon y\right)  }{\left\vert \nabla
_{\mathrm{g}}\theta\left(  \varepsilon y\right)  \right\vert ^{n}}+O\left(
\varepsilon^{1-n}\right)  \label{23}%
\end{equation}
By Lagrange's theorem we can write%
\begin{equation}
\theta\left(  \varepsilon y,\omega\right)  -\theta\left(  0,\omega\right)
=\varepsilon\rho_{\varepsilon}\left(  y\right)  ,\rho_{\varepsilon}\left(
y\right)  =\int_{0}^{1}\left\langle y,\nabla\theta\left(  \varepsilon
ty\right)  \right\rangle \mathrm{d}t \label{16}%
\end{equation}
This yields
\[
\int\log\left\vert \theta\left(  x,\omega\right)  -\theta\left(  x_{0}%
,\omega\right)  \right\vert d_{n}\left(  x\right)  \mathrm{d}x=\int\left(
\log\varepsilon\right)  d_{n}\left(  x\right)  \mathrm{d}x+\int\log\left\vert
\rho_{\varepsilon}\left(  y\right)  \right\vert d_{n}\left(  x\right)
\mathrm{d}x
\]
The first integral vanishes since $d_{n}$ is equal to $x_{1}$-derivatives of a
function with compact support. It follows that the left hand side of
(\ref{22}) equals%
\[
\frac{1}{\left(  n-1\right)  !}\int\log\left\vert \rho_{\varepsilon}\left(
y\right)  \right\vert d_{n}\left(  x\right)  \mathrm{d}x=\int\log\left\vert
\rho_{\varepsilon}\left(  y\right)  \right\vert \frac{1}{\left\vert
\nabla_{\mathrm{g}}\theta\left(  x_{0},\omega\right)  \right\vert ^{n}}%
\frac{\partial^{n}e}{\partial y_{1}^{n}}\mathrm{d}y+O\left(  \varepsilon
\right)
\]
since the logarithmic factor is absolutely integrable. By (\ref{16}) we
have$\ C^{1}$-convergence $\rho_{\varepsilon}\rightarrow\left\langle
y,\nabla\theta\left(  0\right)  \right\rangle $ as $\varepsilon\rightarrow0$
in a neighborhood of the origin. This implies the inequality%
\[
-\int_{\left\vert y\right\vert \leq1,\left\vert y_{1}\right\vert <\delta}%
\log\left\vert y_{1}\right\vert \mathrm{d}y-\int_{\left\vert y\right\vert
\leq1,\rho_{\varepsilon}\left(  y\right)  <\delta}\log\left\vert
\rho_{\varepsilon}\left(  y\right)  \right\vert \mathrm{d}y\leq C\delta
\log\left\vert \delta\right\vert
\]
where $\delta,\ 0<\delta\leq1$ is arbitrary and $C$ does not depend on
$\varepsilon$ and $\delta.$ On the other hand, $\log\left\vert \rho
_{\varepsilon}\left(  y\right)  \right\vert \rightarrow\log\left\vert
y_{1}\right\vert $ everywhere as $\varepsilon\rightarrow0$. Therefore
\[
\int\log\left\vert \rho_{\varepsilon}\left(  y\right)  \right\vert
\frac{\partial^{n}e}{\partial y_{1}^{n}}\mathrm{d}y\rightarrow\int
\log\left\vert y_{1}\right\vert \frac{\partial^{n}e}{\partial y_{1}^{n}%
}\mathrm{d}y
\]
and%
\[
\operatorname{Re}\int\frac{e_{\varepsilon}\left(  x-x_{0}\right)  \mathrm{d}%
V}{\left(  \theta\left(  x,\omega\right)  -\theta\left(  x_{0},\omega\right)
-i0\right)  ^{n}}\rightarrow\frac{1}{\left(  n-1\right)  !\left\vert
\nabla_{\mathrm{g}}\theta\left(  x_{0}\right)  \right\vert ^{n}}\int
\log\left\vert y_{1}\right\vert \frac{\partial^{n}e}{\partial y_{1}^{n}%
}\mathrm{d}y
\]
More detailed arguments show that the difference is equal to $O\left(
\varepsilon^{1/2}\log\varepsilon\right)  .$ The same is true for the linear
function $\theta\left(  x,\omega\right)  =x_{1}$ that is%
\[
\int_{X}\frac{e_{\varepsilon}\left(  x-x_{0}\right)  \mathrm{d}V}{\left(
\theta\left(  x,\omega\right)  -\theta\left(  x_{0},\omega\right)  -i0\right)
^{n}}\rightarrow\frac{1}{\left\vert \nabla_{\mathrm{g}}\theta\left(
x_{0},\omega\right)  \right\vert ^{n}}\operatorname{Re}\int_{\mathbb{R}^{n}%
}\frac{e\left(  y\right)  \mathrm{d}y}{\left(  y_{1}-i0\right)  ^{n}%
},\ \varepsilon\rightarrow0
\]
Calculate the integral in the right hand side by partial integration:%

\begin{align}
\operatorname{Re}\int_{\mathbb{R}^{n}}\frac{e\mathrm{d}y}{\left(
y_{1}-i0\right)  ^{n}}  &  =\frac{1}{n-1}\operatorname{Re}\int_{\mathbb{R}%
^{n}}\frac{\mathrm{d}y}{\left(  y_{1}-i0\right)  ^{n-1}}\frac{\partial
e}{\partial y_{1}}\nonumber\\
&  =\frac{1}{n-1}\operatorname{Re}\int_{\mathrm{S}^{n-1}}\left(  \cos
\omega_{1}-i0\right)  ^{2-n}\mathrm{d}\omega\int_{0}^{\infty}\frac{\partial
e_{0}}{\partial r^{2}}\mathrm{d}r^{2}\nonumber\\
&  =-\frac{1}{n-1}\operatorname{Re}\int_{\mathrm{S}^{n-1}}\left(  \cos
\omega_{1}-i0\right)  ^{2-n}\mathrm{d}\omega\nonumber\\
&  =\frac{\left\vert \mathrm{S}^{n-2}\right\vert }{n-1}\operatorname{Re}%
\int_{0}^{\pi}\sin^{n-2}\omega_{1}\left(  \cos\omega_{1}-i0\right)
^{2-n}\mathrm{d}\omega_{1} \label{18}%
\end{align}
where $y=r\cos\omega_{1}$ since
\[
\frac{\partial e}{\partial y_{1}}=2y_{1}\frac{\partial e}{\partial r^{2}%
},\ y_{1}-i0=\left(  \cos\omega_{1}-i0\right)  r,\ \int_{0}\partial
e_{0}/\partial r^{2}\mathrm{d}r^{2}=-e\left(  0\right)  =-1
\]
By substituting $s=\cos^{2}\omega_{1}$ we get%
\[
\operatorname{Re}\int_{0}^{\pi}\left(  \frac{\sin\omega_{1}}{\cos\omega
_{1}+i0}\right)  ^{n-2}\mathrm{d}\omega_{1}=\frac{1}{2}B\left(  \frac{n-1}%
{2},\frac{3-n}{2}\right)  =\left(  -1\right)  ^{n/2-1}\frac{\pi}{2}%
\]
We use the formula for Beta-function extended for all complex (non negative
integers) values of arguments. The exponent $\lambda=1/2-n/2$ is a regular
point and we can use a classical formula. The right hand side of (\ref{18})
equals%
\[
\left(  -1\right)  ^{n/2-1}\frac{\pi\left\vert \mathrm{S}^{n-2}\right\vert
}{2\left(  n-1\right)  }=\left(  -1\right)  ^{n/2-1}\frac{\pi^{\left(
n+1\right)  /2}}{\Gamma\left(  \left(  n+1\right)  /2\right)  }.
\]
which yields (\ref{19}).

In the case of odd $n$ we integrate by parts as in the previous case and get%
\begin{equation}
\operatorname{Im}\int\frac{e_{\varepsilon}\left(  x-x_{0}\right)  \mathrm{d}%
V}{\left(  \theta\left(  x,\omega\right)  -\theta\left(  x_{0},\omega\right)
-i0\right)  ^{n}}=\frac{\pi}{\left(  n-1\right)  !}\int_{\theta\left(
x,\omega\right)  >\theta\left(  x_{0},\omega\right)  }d_{n}\left(  x\right)
\mathrm{d}x\nonumber
\end{equation}
Taking in account (\ref{23}) and convergence $\rho_{\varepsilon}\rightarrow
y_{1}$ the limit of the right hand side is%
\begin{equation}
\frac{\pi}{\left(  n-1\right)  !\left\vert \nabla_{\mathrm{g}}\theta\left(
x_{0},\omega\right)  \right\vert ^{n}}\int_{y_{1}>0}\frac{\partial^{n}%
e}{\partial y_{1}^{n}}\mathrm{d}y \label{24}%
\end{equation}
Integrating by parts backward gives the equation%
\[
\int_{y_{1}>0}\frac{\partial^{n}e}{\partial y_{1}^{n}}\mathrm{d}y=-\int
_{y_{1}=0}\frac{\partial^{n-1}e}{\partial y_{1}^{n-1}}\mathrm{d}y^{\prime
}=-2^{m-1}\left(  n-2\right)  !!\left\vert \mathrm{S}^{n-2}\right\vert
\int_{0}^{\infty}\frac{\partial^{m}e_{0}\left(  s\right)  }{\partial s^{m}%
}s^{m-1}\mathrm{d}s
\]
where $y^{\prime}=\left(  y_{2},...,y_{n}\right)  ,\ m=\left(  n-1\right)
/2,\ s=r^{2}.$ Here, we applied the formula%
\[
\left.  \frac{\partial^{n-1}e\left(  y\right)  }{\partial y_{1}^{n-1}%
}\right\vert _{y_{1}=0}=2^{m}\left(  n-2\right)  !!\frac{\partial^{m}%
e_{0}\left(  s\right)  }{\partial s^{m}}%
\]
Integrating by parts $m-1$ times in the interior integral we get the quantity
\[
\int_{0}^{1}\frac{\partial^{m}e_{0}\left(  s\right)  }{\partial s^{m}}%
s^{m-1}\mathrm{d}s=\left(  -1\right)  ^{m-1}\left(  m-1\right)  !\int_{0}%
^{1}\partial e_{0}/\partial s\mathrm{d}s=\left(  -1\right)  ^{m}\left(
m-1\right)  !
\]
This implies that\ (\ref{24}) equals
\[
\frac{\left(  -1\right)  ^{m}2^{m-1}\pi\left(  m-1\right)  !\left(
n-2\right)  !!\left\vert \mathrm{S}^{n-2}\right\vert }{\left(  n-1\right)
!}\frac{1}{\left\vert \nabla_{\mathrm{g}}\theta\left(  x_{0},\omega\right)
\right\vert ^{n}}=\frac{2\left(  2\pi i\right)  ^{n-1}}{\left(  n-1\right)
!}\frac{1}{\left\vert \mathrm{S}^{n-1}\right\vert \left\vert \nabla
_{\mathrm{g}}\theta\left(  x_{0},\omega\right)  \right\vert ^{n}}%
\]
For odd $n$ we have%
\[
a_{n}\left(  x_{0},\omega\right)  \doteq\frac{\left(  n-1\right)  !}{\pi}%
\lim_{\varepsilon\rightarrow0}\operatorname{Im}\int_{X}\frac{e_{\varepsilon
}\left(  x-x_{0}\right)  \mathrm{d}V}{\left(  \theta\left(  x,\omega\right)
-\theta\left(  x_{0},\omega\right)  -i0\right)  ^{n}}=\frac{2\left(  2\pi
i\right)  ^{n-1}}{\left\vert \mathrm{S}^{n-1}\right\vert \left\vert
\nabla_{\mathrm{g}}\theta\left(  x_{0},\omega\right)  \right\vert ^{n}}%
\]
Integrating over $\mathrm{S}^{n-1}$ we get
\[
a_{n}\left(  x_{0}\right)  =\int a_{n}\left(  x_{0},\omega\right)
\mathrm{d}\omega=\frac{2\left(  2\pi i\right)  ^{n-1}}{\left\vert
\mathrm{S}^{n-1}\right\vert }\int\frac{\mathrm{d}\omega}{\left\vert
\nabla_{\mathrm{g}}\theta\left(  x_{0},\omega\right)  \right\vert ^{n}}%
\]
and (\ref{2}) follows.$\blacktriangleright$

\section{Integrals of rational trigonometric functions}

We focus now on the condition $\Theta\left(  x,y\right)  =0.$ A function of
the form%
\[
t\left(  \varphi\right)  =\sum_{j=0}^{k}a_{j}\cos j\varphi+b_{j}\sin j\varphi
\]
is called trigonometric polynomial of degree $k$ if $a_{k}\neq0$ or $b_{k}%
\neq0.$ Any trigonometric polynomial is $2\pi$-periodic and is well-defined
and holomorphic on the cylinder $\mathbb{C}/2\pi\mathbb{Z}.$ It always has
$2k$ zeros in the cylinder. If the polynomial is real, number of real zeros is even.

\begin{lemma}
\label{T}Let $t\left(  \varphi\right)  $ and $s\left(  \varphi\right)  $ be
real trigonometric polynomials such that $\mathrm{\deg~}s<\mathrm{\deg~}t$ and
all the roots of $t$ are real. Then for $r=s/t\ $and\ arbitrary natural $n$
\begin{equation}
\left(  P\right)  \int_{0}^{2\pi}r^{n}\left(  \varphi\right)  \mathrm{d}%
\varphi\doteq\frac{1}{2}\int_{0}^{2\pi}\left(  r\left(  \varphi\right)
+i0\right)  ^{n}\mathrm{d}\varphi+\frac{1}{2}\int_{0}^{2\pi}\left(  r\left(
\varphi\right)  -i0\right)  ^{n}\mathrm{d}\varphi=0. \label{7}%
\end{equation}

\end{lemma}

\textbf{Proof.} Suppose first that all roots of $t$ are simple. Let
$\alpha_{1}<\alpha_{2}<...<\alpha_{m}$ be all roots of $\partial
t/\partial\varphi$ on the circle$\ \mathbb{R}/2\pi\mathbb{Z}.$ Let
$\varepsilon_{k}=\mathrm{sgn~}\partial t/\partial\varphi$ on the interval
$\left(  \alpha_{k},\alpha_{k+1}\right)  $ for $k=1,...,m,$ where
$\alpha_{m+1}=\alpha_{1} $. The function $r\left(  \zeta\right)  $ is
meromorphic for $\zeta=\varphi+i\tau\in\mathbb{C}/2\pi\mathbb{Z}$ \ and has no
poles in the half-cylinder $\{\tau>0\}$ since of the assumption. Take a
continuous function $\lambda=\tau\left(  \varphi\right)  $ defined on the
circle that vanishes in all points $\alpha_{k}$ and is positive otherwise. We
have for any $k$
\[
\int_{\alpha_{k}}^{\alpha_{k+1}}\left(  r\left(  \varphi\right)
+\varepsilon_{k}i0\right)  ^{n}\mathrm{d}\varphi=\int_{\alpha_{k}}%
^{\alpha_{k+1}}r^{n}\left(  \varphi+i0\right)  \mathrm{d}\varphi.
\]
Take the sum%
\begin{equation}
\sum_{k}\int_{\alpha_{k}}^{\alpha_{k+1}}\left(  r\left(  \varphi\right)
+\varepsilon_{k}i0\right)  ^{n}\mathrm{d}\varphi=\sum\int_{\alpha_{k}}%
^{\alpha_{k+1}}r^{n}\left(  \varphi+i0\right)  \mathrm{d}\varphi=\int
_{0}^{2\pi}r^{n}\left(  \varphi+i0\right)  \mathrm{d}\varphi\label{6}%
\end{equation}
Now we can replace the form $r^{n}\left(  \varphi+i0\right)  \mathrm{d}%
\varphi$ by $r^{n}\left(  \zeta\right)  \mathrm{d}\zeta$ for $\zeta
=\varphi+i\eta$ for an arbitrary $\eta>0$ without changing the integral in the
right-hand side. We have $\left\vert r\left(  \zeta\right)  \right\vert
\rightarrow0$ as $\eta\rightarrow\infty$ hence the right\ hand side of
(\ref{6}) vanishes. The real part of the left hand side is equal to the
left\ hand side of (\ref{7}) which also vanishes.\ 

In the general case, we can approximate the polynomial $t$ with real roots\ by
polynomials $\tilde{t}$ with real simple roots. The equation (\ref{7}) holds
for $\tilde{r}=s/\tilde{t}$ hence it is true for $r=s/t.$ $\blacktriangleright
$

\begin{proposition}
\label{FO}Let $v\in\mathbb{R}^{n}$ and $a\in\mathbb{R}$ be such that
$\left\vert a\right\vert <\left\vert v\right\vert .$ Then for arbitrary even
$n\geq2$
\[
\operatorname{Re}\int\frac{\mathrm{d}\omega}{\left(  \left\langle
\omega,v\right\rangle -a-i0\right)  ^{n}}=0
\]
and for arbitrary odd $n\geq3$%
\[
\operatorname{Im}\int\frac{\mathrm{d}\omega}{\left(  \left\langle
\omega,v\right\rangle -a-i0\right)  ^{n}}=0.
\]

\end{proposition}

\textbf{Proof.\ }We may assume that $\left\vert v\right\vert =1.$ For even $n
$ we have%
\[
\operatorname{Re}\int\frac{\mathrm{d}\omega}{\left(  \left\langle
\omega,v\right\rangle -a-i0\right)  ^{n}}=\operatorname{Re}\int\frac
{\mathrm{d}\omega}{\left(  \cos\varphi-a-i0\right)  ^{n}}%
\]
where $\varphi$ is the spherical distance between $\omega$ and$\ v.$ We have
$\mathrm{d}\omega=\sin^{n-2}\varphi\mathrm{d}\varphi\mathrm{d}\omega^{\prime}$
where \textrm{d}$\omega^{\prime}$ is the area of a unit sphere $\mathrm{S}%
^{n-2}=\left\{  \omega\in\mathrm{S}^{n-1};\ \varphi=\pi/2\right\}  .$
Integrating over $n-2$-spheres $\varphi=\mathrm{const}$ we get
\[
\operatorname{Re}\int\frac{\mathrm{d}\omega}{\left(  \cos\varphi-a-i0\right)
^{n}}=\frac{\left\vert \mathrm{S}^{n-2}\right\vert }{2}\left(  P\right)
\int_{0}^{2\pi}\frac{\sin^{n-2}\varphi\mathrm{d}\varphi}{\left(  \cos
\varphi-a\right)  ^{n}}%
\]
since the integrand is $\pi$-periodic. The right hand side vanishes by Lemma
\ref{T}.

For odd $n$ we have%
\begin{align*}
\operatorname{Im}\int\frac{\mathrm{d}\omega}{\left(  \cos\varphi-a-i0\right)
^{n}}  &  =\frac{\left\vert \mathrm{S}^{n-2}\right\vert }{2i\left(
n-1\right)  !}\left[  \int_{0}^{\pi}\frac{\sin^{n-2}\varphi\mathrm{d}\varphi
}{\left(  \cos\varphi-a-i0\right)  ^{n}}-\frac{\sin^{n-2}\varphi
\mathrm{d}\varphi}{\left(  \cos\varphi-a+i0\right)  ^{n}}\right] \\
&  =\frac{\left\vert \mathrm{S}^{n-2}\right\vert }{2i\left(  n-1\right)
!}\int_{\left\vert \varphi-\alpha\right\vert =\varepsilon}\frac{\sin
^{n-2}\varphi\mathrm{d}\varphi}{\left(  \cos\varphi-a\right)  ^{n}}=\frac
{\pi\left\vert \mathrm{S}^{n-2}\right\vert }{\left(  n-1\right)
!}\mathrm{res}_{\alpha}\frac{\sin^{n-2}\varphi\mathrm{d}\varphi}{\left(
\cos\varphi-a\right)  ^{n}}%
\end{align*}
where $\alpha=\arccos a\in\left[  0,\pi\right]  .$ Changing variable
$\zeta=\cos\varphi$ and omitting the constant coefficient we get the quantity%
\[
\mathrm{res}_{a}\frac{\left(  1-\zeta^{2}\right)  ^{m}\mathrm{d}\zeta}{\left(
\zeta-a\right)  ^{n}}%
\]
where $m=\left(  n-3\right)  /2.$ The residue is equal to zero since the
numerator has order $2m=n-3.\ \blacktriangleright$

\begin{corollary}
\label{L}For any regular resolved generating function $\Phi=\theta-\lambda$
such that$\ $for\ any\ pair\ of\ points$\ x\neq y$ in an open set$\ \Omega
\subset X\ $we have$\ \theta\left(  x,\omega\right)  -\theta\left(
y,\omega\right)  =\left\langle v,\omega\right\rangle +a,\ \left\vert
a\right\vert <\left\vert v\right\vert ,$ formulas (\ref{9}) and (\ref{2}) hold
for any function $f\in L_{2}\left(  X\right)  $ with support in $\Omega$.
\end{corollary}

\section{Reconstruction in spaces of constant curvature}

We apply the above results to recover few known and unknown inversion formulas
for geodesic integral transforms in spaces of curvature $\kappa=0,1,-1.$

\textbf{Euclidean space. }Take the generating function $\Phi\left(
x;\lambda,\omega\right)  =\left\langle \omega,x\right\rangle -\lambda$ in
$\mathbb{R}^{n}\times\Sigma,~\Sigma=\mathbb{R}\times\mathrm{S}^{n-1}$. We
have$\ \left\vert \nabla\theta\right\vert =D_{n}\left(  x\right)  =1.$ Then
(\ref{9}),(\ref{2}) coincide with the classical reconstruction in Euclidean
space from data of hyperplane integrals.

\textbf{Elliptic space. }Funk\textbf{\ \cite{Fu}} inspired by the seminal
paper of Minkowski \cite{M} found a reconstruction formula of an even function
$f$ on the unit sphere $\mathrm{S}^{2}$ from integrals $Gf$ for the family of
big circles. His method was adapted by Radon for Euclidean plane. A
generalization of Funk's formula for arbitrary dimension is due to Helgason
\cite{He1}.

A reconstruction formula of a different form can be obtained if we apply
Theorem \ref{NR} to a generating function $\Phi\left(  x;\lambda
,\omega\right)  =\left\langle \omega,\xi\right\rangle -\lambda$ defined in
$X\times\Sigma$ where $X=\{\left(  x_{0},x\right)  \in E^{n+1},x_{0}%
^{2}+\left\vert x\right\vert ^{2}=1,x_{0}\geq0\},\ \Sigma=\left\{  \xi\in
E^{n}:\left\vert \xi\right\vert <1\right\}  $ and $\mathrm{g}$ is the metric
in $X$ induced form the Euclidean space $E^{n+1}.\ $Omitting some simple
calculations, we come to

\begin{theorem}
If $n$ is even then any function $f\in L_{2}\left(  X\right)  _{\mathrm{comp}%
}$ can be reconstructed from its integrals $Gf\left(  \sigma\right)  $ over
big spheres $S\left(  \sigma\right)  =\left\{  x\in X,\left\langle
\sigma,x\right\rangle =0\right\}  ,\ \sigma\in\mathrm{S}_{+}^{n}$ by%
\[
f\left(  x\right)  =-\frac{\left(  n-1\right)  !}{\left(  2\pi i\right)  ^{n}%
}\int_{\mathrm{S}_{+}}\frac{Gf\left(  \sigma\right)  \mathrm{d}\sigma
}{\left\langle \sigma,x\right\rangle ^{n}}%
\]
where $\mathrm{S}_{+}^{n-1}=\left\{  \sigma\in\mathbb{R}^{n+1};\left\vert
\sigma\right\vert =1,\sigma_{0}\geq0\right\}  $ is a hemisphere. If $n$ is odd
we have
\end{theorem}

\[
f\left(  x\right)  =\frac{1}{2\left(  2\pi i\right)  ^{n-1}}\int
_{\mathrm{S}_{+}^{n-1}}\delta^{\left(  n-1\right)  }\left(  \left\langle
\sigma,x\right\rangle \right)  Gf\left(  \sigma\right)  \mathrm{d}\sigma.
\]

\textbf{Hyperbolic space\ }$\left(  \kappa=-1\right)  .$ Take the
generating$\ $function $\Phi\left(  x;\lambda,\omega\right)  =\theta
-\lambda,\ \theta=-2\left(  \left\vert x\right\vert ^{2}+1\right)
^{-1}\left\langle \omega,x\right\rangle ,\ -1<\lambda<1\ $in the unit ball
$X\subset\mathbb{R}^{n}$. The hypersurfaces $Z\left(  \lambda,\omega\right)  $
are fully geodesics for the hyperbolic metric $\mathrm{d}_{\mathrm{g}%
}s=2\left(  1-\left\vert x\right\vert ^{2}\right)  ^{-1}\mathrm{d}s$. By a
similar rearrangement we obtain%
\begin{equation}
f\left(  x\right)  =-\frac{\left(  n-1\right)  !}{\left(  2\pi i\right)  ^{n}%
}\int_{Q_{+}}\frac{Gf\left(  \sigma\right)  \mathrm{d}\sigma}{\left\langle
\sigma,x\right\rangle ^{n}} \label{26}%
\end{equation}
for even $n$ and%
\begin{equation}
f\left(  x\right)  =\frac{1}{2\left(  2\pi i\right)  ^{n-1}}\int_{Q_{+}}%
\delta^{\left(  n-1\right)  }\left(  \left\langle \sigma,x\right\rangle
\right)  Gf\left(  \sigma\right)  \mathrm{d}\sigma\label{27}%
\end{equation}
for odd $n$ where $Q_{+}=\left\{  \sigma\in\mathbb{R}^{n+1};\sigma_{0}%
^{2}-\sigma_{1}^{2}-...-\sigma_{n}^{2}=-1,\sigma_{0}\geq0\right\}  $ is the
dual one-fold hyperboloid.

A reconstruction in a different form was done first by Radon \cite{R} $\left(
n=2\right)  $ and Helgason \cite{He1} ($n>2$); formulas (\ref{26}) and
(\ref{27}) are due to Gelfand-Graev-Vilenkin \cite{GGV}.\ 

\section{Equidistant spheres and horospheres in hyperbolic space}

\textbf{Equidistant spheres. }Let $X$ be again a unit $n$ dimensional ball,
$n\geq2\ $and
\[
\Phi\left(  x;\lambda,\omega\right)  =\theta-\lambda,\ \theta\left(
x,\omega\right)  =\frac{p-\left\langle \omega,x\right\rangle }{1-\left\vert
x\right\vert ^{2}},\ \omega\in\mathrm{S}^{n-1}%
\]
be a generating function where $0\leq p<1.$ For a fixed $\omega$ and an
arbitrary $\lambda\neq0$ the hypersurface $Z\left(  \lambda,\omega\right)
=\left\{  x;\ \Phi\left(  x;\lambda,\omega\right)  =0\right\}  $ is the
intersection of $X$ and of an $n-1$ sphere $S\left(  \lambda\right)  $ whereas
$S\left(  0\right)  =\left\{  \left\langle \omega,x\right\rangle -p=0\right\}
$ is a hyperplane; all the spheres $S\left(  \lambda\right)  $ contain $\ n-2$
sphere $S\left(  0\right)  \cap\partial X.$ For arbitrary real $\lambda$ and
$\mu$ the hypersurfaces $S\left(  \lambda\right)  \cap X,~S\left(  \mu\right)
\cap X$ are equidistant with respect to the hyperbolic metric.$\ $Check that
$\Phi$ fulfils the conditions of Corollary \ref{L} for arbitrary $p,~0\leq
p\leq1$. A proof of regularity is a routine. Further we have$\ $%
\[
\theta\left(  x,\omega\right)  -\theta\left(  y,\omega\right)  =-\left\langle
\omega,\frac{x}{1-\left\vert x\right\vert ^{2}}-\frac{y}{1-\left\vert
y\right\vert ^{2}}\right\rangle +p\left(  \frac{1}{1-\left\vert x\right\vert
^{2}}-\frac{1}{1-\left\vert y\right\vert ^{2}}\right)  .
\]
and we need to prove that%
\begin{equation}
\left\vert \frac{x}{1-\left\vert x\right\vert ^{2}}-\frac{y}{1-\left\vert
y\right\vert ^{2}}\right\vert >\left\vert \frac{1}{1-\left\vert x\right\vert
^{2}}-\frac{1}{1-\left\vert y\right\vert ^{2}}\right\vert \label{10}%
\end{equation}
for arbitrary $x\neq y\in X.$\ Squaring both sides we reduce (\ref{10}) to the
obvious inequality $2\left(  1-xy\right)  >\left(  2-\left\vert x\right\vert
^{2}-\left\vert y\right\vert ^{2}\right)  $. The proof is complete and
reconstructions (\ref{9}) and (\ref{2}) follow. Write these formulas in terms
of the geodesic integral transform%
\[
Hf\left(  \sigma\right)  =\int_{Z\left(  \sigma\right)  }f\mathrm{d}%
_{\mathrm{g}}S,\ \sigma=\left(  \lambda,\omega\right)  .
\]
where$\ \mathrm{d}_{\mathrm{g}}S$ is the hyperbolic\ hypersurface element. The
Minkowski-Funk operator $M_{\Phi}$ can be written in terms of the Euclidean
integral transform $Gf\left(  \sigma\right)  =\int_{Z\left(  \sigma\right)
}f\mathrm{d}_{\mathrm{e}}S\ $since of the factorization $\left\vert
\nabla\theta\left(  x,\omega\right)  \right\vert =\left(  1-\left\vert
x\right\vert ^{2}\right)  ^{-1}\sqrt{4\lambda^{2}-4p\lambda+1}$ (see
(\ref{1})). On the other hand,$\ $%
\[
\mathrm{d}_{\mathrm{g}}S=\left(  \frac{2}{1-\left\vert x\right\vert ^{2}%
}\right)  ^{n-1}\mathrm{d}_{\mathrm{e}}S
\]
which yields%
\[
Mf\left(  \lambda,\omega\right)  =\frac{Gf_{1}\left(  \lambda,\omega\right)
}{\sqrt{\lambda^{2}-p\lambda+1/4}}=\frac{Hf_{2}\left(  \lambda,\omega\right)
}{\sqrt{\lambda^{2}-p\lambda+1/4}}%
\]
where%
\[
\ f_{1}\left(  x\right)  =\frac{\left(  1-\left\vert x\right\vert ^{2}\right)
}{2}f\left(  x\right)  ,\ f_{2}\left(  x\right)  =\left(  \frac{1-\left\vert
x\right\vert ^{2}}{2}\right)  ^{n}f\left(  x\right)  .
\]

\begin{corollary}
For any function $f$ with compact support in the unit ball, a reconstruction
is given for even $n$ by%
\begin{equation}
f\left(  x\right)  =-\frac{1}{\left(  4\pi i\right)  ^{n}}\frac{\left(
1-\left\vert x\right\vert ^{2}\right)  ^{2n}}{D_{n}\left(  x\right)  }%
\int_{\mathrm{S}^{n-1}}\int_{\mathbb{R}}\frac{Hf\left(  \lambda,\omega\right)
\mathrm{d}\lambda}{\sqrt{\lambda^{2}-p\lambda+1/4}}\frac{\mathrm{d}\omega
}{\left(  \left\langle \omega,x\right\rangle -\lambda\left(  1-\left\vert
x\right\vert ^{2}\right)  \right)  ^{n}} \label{31}%
\end{equation}
and for odd $n$ by%
\begin{equation}
f\left(  x\right)  =\frac{1}{2\left(  4\pi i\right)  ^{n}}\frac{\left(
1-\left\vert x\right\vert ^{2}\right)  ^{2n}}{D_{n}\left(  x\right)  }%
\int_{\mathrm{S}^{n-1}}\left.  \frac{\partial^{n-1}}{\partial\lambda^{n-1}%
}\frac{Hf\left(  \lambda,\omega\right)  }{\sqrt{\lambda^{2}-p\lambda+1/4}%
}\right\vert _{\lambda=\theta\left(  x,\omega\right)  }\mathrm{d}\omega.
\label{32}%
\end{equation}

\end{corollary}

\textbf{Horospheres. }Taking $p=1$ in the above formulas we obtain the
function
\[
\theta\left(  x,\omega\right)  =\frac{1-\left\langle \omega,x\right\rangle
}{1-\left\vert x\right\vert ^{2}}%
\]
which defines the family of horospheres $\theta\left(  x,\omega\right)
=\lambda,1/2<\lambda<\infty.$ The formulas (\ref{31}) and (\ref{32}) holds for
horospheres if we substitute $\sqrt{4\lambda^{2}-4\lambda+1}=2\lambda-1$ and
integrate in (\ref{31}) over the ray $\left(  1/2,\infty\right)  $ in the
interior integral.

For the cases $n=2$ and $n=3$ reconstruction formulas (of different form) are
contained in the book of Gelfand-Gindikin-Graev \cite{GGG}.

\section{Isofocal hyperboloids}

\textbf{Hyperboloids}. The equation $\lambda=\left\vert x\right\vert
+\varepsilon x_{1},\varepsilon>1$ defines a fold of a two-fold hyperboloid
\[
\left(  \alpha x_{1}-\frac{\varepsilon\lambda}{\alpha}\right)  ^{2}-x_{2}%
^{2}-...-x_{n}^{2}=\frac{\lambda^{2}}{\alpha^{2}},\ \alpha=\sqrt
{\varepsilon^{2}-1}%
\]
with a focus at the origin. The function $\Phi\left(  x;\lambda,\omega\right)
=\theta\left(  x,\omega\right)  -\lambda,$ $\theta\left(  x,\varphi\right)
=\left\vert x\right\vert +\varepsilon\left\langle x,\omega\right\rangle $
generates the family of all one-fold hyperboloids in $X=\mathbb{R}%
^{n}\smallsetminus0$ with a focus at the origin. The function%
\[
\theta\left(  x,\omega\right)  -\theta\left(  y,\omega\right)  =\varepsilon
\left\langle \omega,x-y\right\rangle +\left\vert x\right\vert -\left\vert
y\right\vert ,\ \omega\in\mathrm{S}^{n-1}%
\]
satisfies the conditions of Proposition \ref{FO} since $\left\vert \left\vert
x\right\vert -\left\vert y\right\vert \right\vert <\varepsilon\left\vert
x-y\right\vert $ for any $x,y\in\mathbb{R}^{n},$ $x\neq y.$ Therefore Theorem
\ref{NR} holds for this family. We have%
\[
\left\vert \nabla\theta\left(  x,\omega\right)  \right\vert ^{2}%
=1+\varepsilon^{2}+2\varepsilon\left\vert x\right\vert ^{-1}\left\langle
\omega,x\right\rangle
\]
and%
\[
D_{n}=\int_{\mathrm{S}^{n-1}}\frac{\mathrm{d}\omega}{\left(  1+\varepsilon
^{2}+2\varepsilon\left\vert x\right\vert ^{-1}\left\langle \omega
,x\right\rangle \right)  ^{n/2}}=\left\vert \mathrm{S}^{n-2}\right\vert
\int_{0}^{\pi}\frac{\sin^{n-2}\varphi\mathrm{d}\varphi}{\left(  1+\varepsilon
^{2}-2\varepsilon\cos\varphi\right)  ^{n/2}}%
\]
where $\varphi$ is the angle between $\omega$ and $\left\vert x\right\vert
^{-1}x$. This integral does not depend on $x.$

\begin{corollary}
For any smooth function $f$ with compact support in $\mathbb{R}^{n}$ and even
$n,$ the equation holds%
\begin{align*}
f\left(  x\right)   &  =-\frac{1}{\left(  2\pi i\right)  ^{n}D_{n}}%
\int_{\mathrm{S}^{n-1}}\int_{\mathbb{R}}\frac{Mf\left(  \lambda,\omega\right)
\mathrm{d}\lambda\mathrm{d}\omega}{\left(  \left\vert x\right\vert
+\varepsilon\left\langle \omega,x\right\rangle -\lambda\right)  ^{n}%
},\ \ n\text{ even}\\
f\left(  x\right)   &  =\frac{1}{2\left(  2\pi i\right)  ^{n-1}D_{n}}%
\int_{\mathrm{S}^{n-1}}\frac{\partial^{n-1}}{\partial\lambda^{n-1}}\left.
Mf\left(  \lambda,\omega\right)  \right\vert _{\lambda=\left\vert x\right\vert
+\varepsilon\left\langle x,\omega\right\rangle }\mathrm{d}\omega
,\ \ n\ \ \text{odd.}%
\end{align*}

\end{corollary}

\section{Photoacoustic geometries}

Consider a resolved generating function%

\[
\Phi\left(  x;\lambda,\omega\right)  =\left\vert x-\xi\left(  \omega\right)
\right\vert ^{2}-\lambda,\omega\in\mathrm{S}^{n-1}%
\]
where $\xi:\mathrm{S}^{n-1}\rightarrow\mathbb{R}^{n}$ is a smooth map. We call
the image $\mathbf{C}$ of $\xi$ central set. Any surface $Z\left(
\lambda,\omega\right)  $ is a sphere of radius $\sqrt{\lambda}$ with the
center $\xi\left(  \omega\right)  \in\mathbf{C}$ and
\[
M_{\Phi}f\left(  \lambda,\omega\right)  =\frac{Gf\left(  \lambda
,\omega\right)  }{2\sqrt{\lambda}}%
\]
where $Gf\left(  \lambda,\omega\right)  \ $is the Euclidean integral over this
sphere. Inversion of the operator $M_{\Phi}$ implies inversion of the
spherical integral transform $G$ for the given central surface $\mathbf{C}$
(and vice versa). This subject is of special interest in view of application
to the photoacoustic (thermoacoustic) tomography, see surveys \cite{IP}%
,\cite{KK}. Inversion formulas for a function supported in a half space with
the hyperplane central set is found by Fawcett \cite{Fa}. In \cite{P0} a
reconstruction was done by reduction to the Radon transform. For a spherical
central surface Finch with coauthors \cite{F1},\cite{F2} found a
reconstruction formula of type (\ref{33})-(\ref{34}) in the physical domain
for arbitrary dimension. Another reconstruction formula was proposed by
Kunyanski \cite{K1}; it is similar to (\ref{33})-(\ref{34}) after a
simplification. An inversion of the spherical mean operator in three
dimensional space was constructed by Xu and Wang \cite{XW} in terms of a
Dirichlet-Green function. This approach gives an explicit reconstruction in
the frequency domain when $\mathbf{C}$ is a sphere or a circular cylinder.
Kunyanski \cite{K2} construct inversion for polyhedral center sets with
special symmetries. Natterer \cite{N} has found an explicit inversion in the
physical domain when the central surface $\mathbf{C}$ is an arbitrary
ellipsoid. Surveys of related results are done in \cite{F2},\cite{XW}%
,\cite{K2}. We show here that for an arbitrary ellipsoid (elliptical cylinder)
in $R^{n}$ a simple reconstruction is given by Theorem \ref{NR}.

\textbf{Ellipsoids.} Set $\xi\left(  \omega\right)  =\left(  a_{1}\omega
_{1},...,a_{n}\omega_{n}\right)  $ where $a_{1},...,a_{n}$ are positive
constants. The central hypersurface $\left\{  x=\xi\left(  \omega\right)
,\omega\in\mathrm{S}^{n-1}\right\}  \ $is the boundary of an ellipsoid $E_{a}$
with half axes $a_{1},...,a_{n}.$ Then
\[
\theta\left(  x,\omega\right)  -\theta\left(  y,\omega\right)  =2\left\langle
\xi\left(  \omega\right)  ,y-x\right\rangle +\left\vert x\right\vert
^{2}-\left\vert y\right\vert ^{2}=2\left\langle \omega,z\right\rangle
+\left\vert x\right\vert ^{2}-\left\vert y\right\vert ^{2}%
\]
where$\ z=\left(  a_{1}\left(  y_{1}-x_{1}\right)  ,...,a_{n}\left(
y_{n}-x_{n}\right)  \right)  .$ The inequality holds
\[
\left\vert \left\vert x\right\vert ^{2}-\left\vert y\right\vert ^{2}%
\right\vert =\left\vert \sum\left(  y_{i}-x_{i}\right)  \left(  y_{i}%
+x_{i}\right)  \right\vert \leq\sum\left\vert a_{i}\left(  y_{i}-x_{i}\right)
\right\vert \sum\left\vert a_{i}^{-1}\left(  y_{i}+x_{i}\right)  \right\vert
\leq\left\Vert z\right\Vert \left\Vert w\right\Vert
\]
where $w=\left(  a_{1}^{-1}\left(  x_{1}+y_{1}\right)  ,...,a_{n}^{-1}\left(
x_{n}+y_{n}\right)  \right)  .$ Suppose that $x,y\in E_{a}$ and $x\neq y$
then\ the point $\left(  x+y\right)  /2$ belongs to the interior of $E_{a}$
which implies $\left\Vert w\right\Vert <2.$ It follows that the right hand
side of\ is \ strictly bounded by $2\left\Vert z\right\Vert $. By Proposition
\ref{FO} Theorem \ref{NR} holds for any $n\geq2$. Note that$\ \left\vert
\nabla\theta\right\vert =2\left\vert x-\xi\left(  \omega\right)  \right\vert
=2\sqrt{\lambda}.$ It follows that any function $f$ supported in the closed
ellipsoid $E_{a}$ can be reconstructed by the formula%
\begin{equation}
f\left(  x\right)  =\frac{1}{\left(  2\pi i\right)  ^{n}D_{n}\left(  x\right)
}\int_{\mathrm{S}^{n-1}}\left(  P\right)  \int_{\mathbb{R}}\frac{Gf\left(
\rho^{2},\omega\right)  \mathrm{d}\rho\mathrm{d}\omega}{\left(  \left\vert
x-\xi\left(  \omega\right)  \right\vert ^{2}-\rho^{2}\right)  ^{n}} \label{33}%
\end{equation}
for even $n$ where we did the substitution $\lambda=\rho^{2},$ and by%
\begin{equation}
f\left(  x\right)  =\frac{1}{4\left(  2\pi i\right)  ^{n-1}D_{n}\left(
x\right)  }\int_{\mathrm{S}^{n-1}}\left(  \frac{\partial}{\partial\rho^{2}%
}\right)  ^{n-1}\left.  \frac{Gf\left(  \rho^{2},\omega\right)  }{\rho
}\right\vert _{\rho=\left\vert x-\xi\left(  \omega\right)  \right\vert
}\mathrm{d}\omega\label{34}%
\end{equation}
for odd $n$ where%

\begin{equation}
D_{n}\left(  x\right)  =\frac{1}{2^{n}\left\vert \mathrm{S}^{n-1}\right\vert
}\int\frac{\mathrm{d}\omega}{\left\vert x-\xi\left(  \omega\right)
\right\vert ^{n}} \label{35}%
\end{equation}

\textbf{Cylinders. }If a central set $\mathbf{C}$ is unbounded can not apply
the same method since it can not be regularly parametrized. However one can
write a reconstruction formula for any closed cylinder $\mathrm{E}$ with an
elliptic base. Indeed, \textrm{E} is a union of the family of ellipsoids
$E_{a}$ as one or several half axes, say $a_{1},...,a_{p},$ tend to infinity,
$a_{p+1},...,a_{n}$ being fixed. One can come to limits in (\ref{33}) and in
(\ref{34}). We omit details.

\textbf{Algebraic plane curves. }In the case $n=2,$ there are more geometries
which admit exact reconstruction formulas. We call a curve $\mathbf{C}%
\subset\mathbb{R}^{2}$ \textit{trigonometric }of degree $k$ if it is given by
a parametric equation%
\[
x_{1}=\xi_{1}\left(  \varphi\right)  ,x_{2}=\xi_{2}\left(  \varphi\right)
,\ \varphi\in\mathrm{S}^{1}%
\]
where $\xi_{1},\xi_{2}$ are real trigonometric polynomials of degree $k$. A
trigonometric curve is always a component of a real algebraic curve. A point
$x\in\mathbb{R}^{2}$ is called \textit{hyperbolic }with respect to a
trigonometric curve $\mathbf{C}$\ if any straight line $L$ through $x$ meets
the curve at $2k$ points (counting with multiplicities). It is easy to see
that the set $H$ (called hyperbolic set) of all hyperbolic points\ is always
closed and convex. Introduce an Euclidean structure in $\mathbb{R}^{2}$ and
consider a function%
\[
\theta\left(  x,\varphi\right)  =\left\vert x-\xi\left(  \varphi\right)
\right\vert ^{2},\ \xi\left(  \varphi\right)  =\left(  \xi_{1}\left(
\varphi\right)  ,\xi_{2}\left(  \varphi\right)  \right)  ,\ \varphi
\in\mathrm{S}^{1},\ x\in\mathbb{R}^{2}.
\]

\begin{proposition}
Let $H$ be the set of hyperbolic points with respect to a trigonometric curve
of degree $k$. For arbitrary points $x,y\in H,\ x\neq y,\ $all roots of the
polynomial $\theta\left(  x,\varphi\right)  -\theta\left(  y,\varphi\right)  $
(of order $k$) are real.
\end{proposition}

\textbf{Proof.} We have%
\begin{align}
\theta\left(  x,\varphi\right)  -\theta\left(  y,\varphi\right)   &
=\left\vert x-\xi\left(  \varphi\right)  \right\vert ^{2}-\left\vert
y-\xi\left(  \varphi\right)  \right\vert ^{2}=\left\vert x\right\vert
^{2}-\left\vert y\right\vert ^{2}-2\left\langle x-y,\xi\left(  \varphi\right)
\right\rangle \nonumber\\
&  =2\left\langle x-y,\xi\left(  \varphi\right)  -s\right\rangle \label{29}%
\end{align}
where $s=\left(  x+y\right)  /2.$ This point is contained in $H$ since $H$ is
convex. Therefore, the line $L=\left\{  z=s+rv,r\in\mathbb{R}\right\}  $ has
$2k$ common points $\xi\left(  \varphi_{1}\right)  ,...,\xi\left(
\varphi_{2k}\right)  $ with $\mathbf{C}$ for arbitrary vector $v\neq0.$ If $v$
is orthogonal to $x-y$ the right hand side of (\ref{29}) vanishes. The
corresponding angles $\varphi_{1},...,\varphi_{2k}$ are real roots of the
polynomial $\theta\left(  x,\varphi\right)  -\theta\left(  y,\omega\right)
$.$\ \blacktriangleright$

The family of circles centered at the curve $\mathbf{C}$ is generated by the
function $\Phi\left(  x;\lambda,\varphi\right)  =\theta\left(  x,\varphi
\right)  -\lambda.\ $Applying Lemma \ref{T} we get

\begin{corollary}
Let $\mathbf{C}\ $be a trigonometric curve with a hyperbolic set $H$. Theorem
\ref{NR} holds for the family of circles centered at $\mathbf{C}$, arbitrary
function $f$ supported in$\ H$ and for both integral transforms $M_{\Phi}$ and
$G.$
\end{corollary}

There is a large variety of trigonometric curves $\mathbf{C}$ with non empty
hyperbolic sets.

\textbf{Examples.\ 1.\ }Let$\ \xi_{1}\left(  \varphi\right)  =2\cos
2\varphi-\cos\varphi,\ \xi_{2}\left(  \varphi\right)  =2\sin2\varphi
+\sin\varphi.$ The curve $\mathbf{C}$ is shown in Fig.1%
\begin{center}
\fbox{\includegraphics[
natheight=3.000900in,
natwidth=4.498700in,
height=3.0009in,
width=4.4987in
]%
{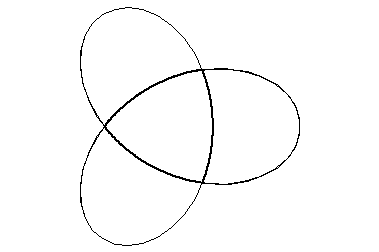}
}\\
Fig.1
\end{center}

The hyperbolic set $H\ $is the triangle in the middle.

\textbf{2. }A hyperbolic "square" set is defined by the trigonometric curve
$\xi_{1}\left(  \varphi\right)  =2\cos3\varphi+\cos\varphi,\ \xi_{2}\left(
\varphi\right)  =2\sin3\varphi-\sin\varphi,$ see Fig.2%

\begin{center}
\fbox{\includegraphics[
natheight=3.000500in,
natwidth=4.499500in,
height=3.0005in,
width=4.4995in
]%
{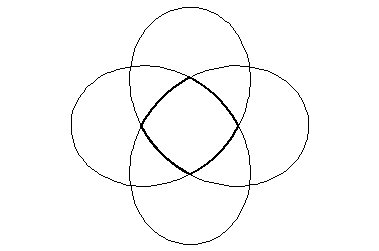}
}\\
Fig.2
\end{center}

\textbf{3.} A "pentagon" is the hyperbolic set of the curve $\xi_{1}\left(
\varphi\right)  =5\cos4\varphi+4\cos\varphi,\ \xi_{2}\left(  \varphi\right)
=5\sin4\varphi-4\sin\varphi,$ see Fig.3:
\begin{center}
\fbox{\includegraphics[
natheight=3.000500in,
natwidth=4.499500in,
height=3.0005in,
width=4.4995in
]%
{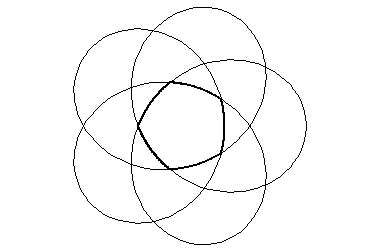}
}\\
Fig.3
\end{center}
etc.

\end{document}